\documentclass[11pt]{article}
\usepackage{amsfonts}
\usepackage{amsmath,amssymb,amsthm,enumerate,epsfig,graphicx}
\usepackage{ifthen,latexsym,syntonly}
\usepackage{natbib}

\rmfamily
\newtheorem{theorem}{Theorem}
\newtheorem{corollary}{Corollary}

\newtheorem{definition}{Definition}

\newtheorem{lemma}{Lemma}
\newtheorem{proposition}{Proposition}
\newtheorem{remark}{Remark}
\newtheorem{assumption}{Assumption}
\newtheorem{condition}{Condition}
\def\baselinestretch{1.0}
\renewcommand{\cite}{\citet*}

\def\G{{\mathcal G}}
\def\XX{\left\langle X,X\right\rangle}
\def\XXhat{\widehat{\left\langle X,X\right\rangle }}

\def\blackbox{\vrule height4pt width3pt depth0pt}
\def\End{\null\nobreak\hfill\blackbox\par\medbreak}
\begin{document}

\title{Efficient Estimation of Stochastic Volatility Using Noisy
  Observations: A Multi-Scale Approach 
\thanks{
Lan Zhang is Assistant Professor, Department of Finance, University
of Illinois at Chicago, Chicago, IL 60607 (E-mail: lanzhang@uic.edu), and
Assistant Professor, Department of Statistics, Carnegie Mellon
University, Pittsburgh, PA 15213.
I gratefully acknowledge the support of the National Science Foundation under
grant DMS-0204639. I would like to thank the Referees and the Editor
for suggestions which greatly improved the paper.
}}

\date{First version: August 15, 2004. \\ This version: December 29, 2005 }

\author{Lan Zhang}
\maketitle

\begin{abstract}
With the availability of high frequency financial data, nonparametric
estimation of volatility of an asset return process becomes
feasible. A major problem is how to estimate the volatility
consistently and efficiently, when the observed asset
returns contain error or noise, for example, in the form of 
microstructure noise. The former (consistency) has been addressed
 in the recent literature. However, the resulting
estimator is not efficient. In \cite{twoscales}, the best estimator
converges to the true volatility only at the rate of $n^{-1/6}$. In this
paper, we propose an estimator, the {\it Multi-scale Realized
Volatility (MSRV)}, which converges to the true volatility
at the rate of $n^{-1/4}$, which is the best attainable. We have shown a central limit theorem for the MSRV estimator, which permits setting intervals for the true integrated volatility on the basis of MSRV.

\medskip

{\it Some key words and phrases}: \
           consistency, dependent noise, discrete observation, efficiency, 
It{\^o} process, microstructure noise, observation error, rate of convergence, 
realized volatility

\end{abstract}

\thispagestyle{empty} \newpage

\markright{Multi Scale Realized Volatility}

\pagestyle{myheadings}

\setcounter{page}{1}

\renewcommand{\baselinestretch}{1.1}
\normalsize%
\section{Introduction}
This paper is about how to estimate volatility non-parametrically and
efficiently. 

With the availability of high frequency financial data, nonparametric
estimation of volatility of an asset return process becomes
feasible. A major problem is how to estimate the volatility
consistently and efficiently, when the observed asset
returns are noisy. The former (consistency) has been addressed
 in the recent literature. However, the resulting
estimator is not  efficient. In \cite{twoscales}, the best estimator
converges to the true volatility only at the rate of $n^{-1/6}$. In this
paper, we propose an estimator which converges to the true
volatility at the rate of $n^{-1/4}$, which is the best attainable. The 
new estimator remains consistent
when the observation noise is dependent.
We call the estimator the {\it Multi Scale Realized Volatility (MSRV)}

To demonstrate the idea, consider $\{Y\}$ as the observed log prices of
a financial instrument, and the
observations take place at the grid of time points 
$\G_n = \{ t_{n,i} , i=0, 1, 2, \cdots n \}$ that span the time interval $[0,T]$.
For the purposes of asymptotics, we shall let $\G_n$ become dense
in $[0,T]$ as $n \rightarrow \infty$. 

Suppose that $\{Y_{t_{n,i}}\}$ are noisy, the corresponding the true (latent) log prices
are $\{X\}$. 
Their relation can be modeled as,
\begin{equation}
Y_{t_{n,i}}=X_{t_{n,i}}+\epsilon_{t_{n,i}}.\label{eq:addModel}
\end{equation}
where $t_{n,i} \in \G_n$. The noise $\epsilon_{t_{n,i}}$s will be assumed to be independent of $X$ and iid. 

The model in (\ref{eq:addModel}) is quite realistic, as evidenced
by the existence of microstructure noise in the price process (\cite{brown90},
\cite{binzhou96}, \cite{corsizumbachmullerdacorogna01}).

We further assume that the true log prices $\{X\}$ satisfy the following equation: 
\begin{equation}
dX_t=\mu_t dt+\sigma_t dB_t \label{eq:diffusion}
\end{equation}
where $B_{t}$ is a standard Brownian motion. Typically, the drift coefficient $\mu_t$ and the diffusion coefficient
$\sigma_t$ are stochastic in the sense that 
\begin{equation}
dX_t (\omega)=\mu (t,\omega) dt+\sigma (t,\omega) dB_t (\omega) \label{eq:Ito}
\end{equation}

Throughout this paper, we use the notation in (\ref{eq:diffusion})
to  denote (\ref{eq:Ito}). By the model in (\ref{eq:Ito}), we mean
that $\{X\}$ follows  an It\^{o} process. A special case is that
$\{X\}$ is Markov, where $\mu_t=\mu(t,X_t)$, and
$\sigma_t=\sigma(t,X_t)$.  In financial literature, $\sigma_t$ is
called the instantaneous {volatility} of $X$.

Our goal is to estimate $\int_{0}^{T}\sigma_{t}^{2}dt$, where
$T$ can be a day, a month, or other time horizon(s). For simplicity, we
call $\int_{0}^{T}\sigma_{t}^{2}dt$ the integrated volatility, and
denote it by
\begin{equation*}
\XX =\int_{0}^{T}\sigma_{t}^{2}dt.
\end{equation*}

The general question is, how to  estimate nonparametrically 
$\int_{0}^{T}\sigma_{t}^{2}dt$, if one can only observe the noisy data
$Y_{t_{n,i}}$ at discrete times  $t_{n,i}\in \G_n$. $\G_n$ is formally defined in Section \ref{sec:as-dis}.

To the best of our knowledge, there are two types of nonparametric estimators for
$\int_{0}^{T}\sigma_{t}^{2}dt$ in the current literature. 
The first type, the simpler one, is  to sum up all  the squared returns in $[0, T]$: 
\begin{equation}
[Y,Y]^{(n,1)} =\sum_{t_{n,i}\in \G_n, i\ge 1}(Y_{t_{n,i}}-Y_{t_{n,{i-1}}
})^{2}, \label{eq:revolall}
\end{equation}
this estimator 
is
generally called {\it realized volatility} or {\it realized
variance} (or {\it RV} for short). However, 
it has been reported that realized volatility
using high-frequency data is not desirable 
(see, for example, \cite{brown90}, \cite{binzhou96},
\cite{corsizumbachmullerdacorogna01}
). 
The reason is that it is not consistent,
even if the noisy observations $Y$ are available continuously. Under 
 discrete observations, 
the bias and the variance of the realized volatility are the same order as the sample size 
$n$.

A slight modification of (\ref{eq:revolall}) is to use the sum of 
squared returns from a  ``sparsely selected'' sample, that is, using 
a subgrid of $\G_n$.
The idea is that by using sparse data, one reduces
the bias and the variance of the conventional realized volatility. This approach is quite popular in the empirical
finance literature. However, this ``sparse'' estimator is
still not consistent, 
 in addition,
 which data to subsample and which
to discard is arbitrary. The behavior of this type of estimator,
and a sufficiency based improvement of it, is analyzed in \cite{twoscales}.

A second type of estimator for $\int_{0}^{T}\sigma_{t}^{2}dt$ is based
on {\it two sampling scales}. As introduced in Section 4 (p. 1402) of
\cite{twoscales}, the Two Scales Realized Volatility {\it (TSRV)} has
the form 
\begin{eqnarray}
\XXhat^{(TSRV)}=[Y,Y]^{(n,K)}-2\frac{n-K+1}{nK}[Y,Y]^{(n,1)}, \label{eq:revol2scale}
\end{eqnarray}
where 
\begin{eqnarray}
\label{eq:revolK1}
[Y,Y]^{(n,K)}&=& \frac{1}{K} \sum_{t_{n,i}\in \G_n, i\ge K} {(Y_{t_{n,i}}
-Y_{t_{n,{i-K}}})}^2,
\end{eqnarray}
with $K$ being a positive integer. 
Thus the estimator in (\ref{eq:revol2scale}) averages 
the squared  returns from sampling every data point ($[Y,Y]_{T}^{(n,1)}$)
and those from  sampling every $K$-th data point ($[Y,Y]_{T}^{(n,K)}$). Its asymptotic behavior
was derived when $K \rightarrow \infty$
as $n \rightarrow \infty$. The TSRV estimator has many desirable features, including asymptotic unbiasedness, consistency, and asymptotic normality\footnote{
A related estimator can be found in \cite{binzhou96} and \cite{hansenlunde06}, however, their estimator 
(takes $k$ to be fixed) does not yield a consistent estimator.}.
However, its rate of convergence
is not satisfactory. For an instance, the best estimator in \cite{twoscales}
converges to $\int_{0}^{T}\sigma_{t}^{2}dt$ at the rate of $n^{-1/6}$.

In this paper, we propose a new class of estimators, collectively
referred to as {\it Multi Scale Realized Volatility (MSRV)}
which converge to
$\int_{0}^{T}\sigma_{t}^{2}dt$ at the rate of  $n^{-1/4}$. This new
estimator has the form,
\begin{equation*}
\XXhat^{(n)} = \sum_{i=1}^{M} \alpha_i [Y, Y]^{(n,K_i)}.
\end{equation*}
where $M$ is a positive integer greater than 2. Comparing to $\widehat{ \XX}_{T}^{(TSRV)}$ which uses two time scales (1
and $K$), $\XXhat^{(n)}$ combines $M$ different time scales. The weights
$a_i$ are selected so that $\XXhat^{(n)}$ is unbiased and has optimal
convergence rate. The rationale is that by combining more than two time
scales, we can improve the efficiency of the estimator. Interestingly,
the $n^{-1/4}$ rate of convergence in our new estimator is the same as
the one in parametric estimation for volatility, when the true process
is Markov (see \cite{gloterjacod} ). Thus this rate is the best
attainable. 
Earlier related results in the same direction can be found
in \cite{stein87, stein90, stein93} and \cite{ying91, ying93}. See
also \cite{yacpmrfs05}. Related independent work can also be found in 
\cite{barndorffnielsenhansenlundeshephard}. For the estimating functions-based approach,
there is a nice review by  \cite{sorensenhandbook}.

We emphasize that our MSRV estimator is
nonparametric, and the true process follows a more general It{\^o}
process, where the volatility could depend on the entire history of
the $X$ process plus additional randomness. 

The paper is organized as following. In section 2, we motivate the
idea of averaging over $M$ different time scales. As we shall see, our
estimator is unbiased, and its asymptotic variance comes from  the
noise (the $\epsilon_{t_{n,i}}$s) as well as from  the discreteness of the sampling
times $t_{n,i}$.
In Sections 3-4, we derive the weights $a_i$'s which are optimal for
minimizing the variance that comes from 
noise, and we give a central limit theorem for the contribution of the
noise term. A specific family of weights is introduced in section 4.
We then elaborate on the discretization error in Section 5,
and show a CLT for this error. 
Section 6 the gives the central limit theorem for the MSRV estimator.

For the statements of results, we shall use the following assumptions:
\begin{assumption}
\label{ass:process}
(Structure of the latent process).
The $X$ process is adapted to a filtration $({\cal X}_t)$, and
satisfies
(\ref{eq:diffusion}), where $B_t$ is an $({\cal X}_t)$-Brownian motion, and the
$\mu_t$ and $\sigma_t$ are 
$({\cal X}_t)$-adapted processes which are continuous almost surely.
Also both processes are bounded above by a constant, and
$\sigma_t$ is bounded away from zero. We denote ${\cal X} = {\cal X}_T$.

As a technical matter,
we suppose that there is a $\sigma$-field ${\cal N}$
and a continuous finite dimensional
local martingale $(M_t)$ so that 
${\cal X}_t = \sigma(M_s , 0 \le s \le t) \vee {\cal N}$.
\end{assumption}

\begin{assumption}
\label{ass:noise}
(Structure of the noise). The $\epsilon_{t_{n,i}}$ are independent and identically distributed,
with $E[\epsilon] = 0$ and $E[\epsilon^4] < \infty$. The $\epsilon_{t_{n,i}}$ are 
also independent of ${\cal X}$
\end{assumption}

These assumptions are not minimal for all results. In terms
of the structure of the process,
see, for example,
Section 5 in \cite{jacodprotter98} and
Proposition 1 in \cite{myklandzhang} for examples of statements 
where the $\mu$ and $\sigma$ processes are not assumed to be continuous.
For the methodology to incorporate dependence into the noise structure,
see \cite{yacperlants}. Our current assumptions, however, provide
a setup with substantial generality without overly complicating
the proofs. 

The final item in Assumption \ref{ass:process} is standard for
the type of limit result that we discuss, cf. similar conditions
in \cite{jacodprotter98}, \cite{zhang2001}, \cite{myklandzhang} and
\cite{twoscales}.

\section{Motivation: Averaging the Observations
of $\XX$}

In \cite{twoscales}, we have observed that by combining the square increments of the
returns from two time scales, the resulting two-scale estimator
$\widehat{\XX_T}^{(TSRV)}$ in (\ref{eq:revol2scale}) improves upon the
realized volatility, which uses only one time scale, as in
(\ref{eq:revolall}). The improvement is about
reducing both the bias and the variance.

If the two-scale estimator is better than the one-scale estimator, a
natural question would be how about the estimator combining more than
2 time scales. This question motivates the present paper. In this
section we briefly
go through the main argument. 
  
To proceed, recall definition (\ref{eq:revolK1}) of $[Y,Y]^{(n,K)}$,
and set, similarly,
\begin{equation}
\label{eq:revolK2}
[X,\epsilon]^{(n,K)} = \frac{1}{K} \sum_{t_{n,i}\in \G_n, i\ge K} (X_{t_{n,i}}
-X_{t_{n,{i-K}}})(\epsilon_{t_{n,i}} -\epsilon_{t_{n,{i-K}}}) ,
\end{equation}
and
\begin{eqnarray*}
[\epsilon,\epsilon]^{(n,K)} = \frac{1}{K} \sum_{t_{n,i}\in \G_n, i\ge K} {(\epsilon_{t_{n,i}} -\epsilon_{t_{n,{i-K}}})}^2 .
\end{eqnarray*}
Under (\ref{eq:addModel}), one can decompose $[Y, Y]^{(n,K)}$ into 
$$[Y, Y]^{(n,K)}=[X, X]^{(n,K)}+ [\epsilon, \epsilon]^{(n,K)}+2[X,\epsilon]^{(n,K)}
.$$

We consider estimators on the form
\begin{equation}
\XXhat^{(n)} = \sum_{i=1}^{M} \alpha_i [Y, Y]^{(n, K_i)} 
\end{equation}
where $\alpha_i$'s are the weights to
be determined. A first intuitive requirement is obtained by noting that
\begin{equation}
E( \XXhat^{(n)} | X \mbox{ process } ) = 
\sum_{i=1}^{M} \alpha_i  [X, X]^{(n, K_i)} 
+ 2 E \epsilon^2 \sum_{i=1}^{M} \alpha_i \frac{n+1-K_i}{K_i}
\end{equation}
Since $[X, X]^{(n,K_i)}$ are asymptotically unbiased for $\XX$
(\cite{twoscales}), it is natural to require that
\begin{equation}
\sum_{i=1}^{M} \alpha_i = 1 \mbox{ and } \sum_{i=1}^{M} \alpha_i \frac{n+1-K_i}{K_i} =0
\end{equation}
A slight redefinition will now make the problem more transparent.
Let
\begin{equation}
a_1 = \alpha_1 - \left [ (n+1) \left(\frac{1}{K_1} - \frac{1}{K_2} \right ) \right ]^{-1} , ~~
a_2 = \alpha_2 - (a_1 - \alpha_1 ) \mbox{ and } a_i =
\alpha_i \mbox{ for } i \ge 3.
\end{equation}
Our conditions on the $\alpha$'s are now equivalent to
\begin{condition}\label{cond:ace1}
$\sum a_i =1,$
\end{condition}
\begin{condition}\label{cond:ace2}
$\sum_{i=1}^{M} \frac{a_i}{K_i} =0.$
\end{condition}
\noindent
To understand the estimator $\XXhat^{(n)}$ in terms of the $a_i$'s, consider the
following asymptotic statement. Here, and everywhere below, we
allow $a_i$, $K_i$ and $M$ to depend on $n$ ({\it i.e.}, they have the form
$a_{n,i}$, $K_{n,i}$ and $M_n$), though sometimes the dependence
on $n$ is suppressed in the notation. We obtain (for proof, see Section \ref{sec:proofs})
\begin{proposition}
\label{prop:approxA}
Suppose that $K_{n,1}$ and $K_{n,2}$ are $O(1)$ as $n \rightarrow \infty$.
Under Assumptions \ref{ass:process}-\ref{ass:noise},
\begin{equation}
\XXhat^{(n)} = \sum_{i=1}^{M} a_{i} [Y, Y]^{(n,K_{i})} - 2 E \epsilon^2 +
O_p(n^{-1/2})
\label{temp1a}
\end{equation}
\end{proposition}

To further analyze the terms in (\ref{temp1a}), write
\begin{equation}
\label{eq:approx1}
[Y, Y]^{(n,K)} = [X, X]^{(n,K)} + \frac{2}{K} \sum_{i=0}^{n} \epsilon_{t_{n,i}}^2
+U_{n,K}+ V_{n,K}
\end{equation}
where $U_{n,K}$ will turn out to be the main error term,
\begin{equation}
U_{n,K} = - \frac{2}{K} \sum_{i=K}^{n} \epsilon_{t_{n,i}} \epsilon_{t_{n,{i-K}}},
\end{equation}
and $V_{n,K}$ will be
a remainder term, given by
$V_{n,K}= 2 [X, \epsilon]^{(n,K)}  -\frac{1}{K}\sum_{i=0}^{K-1}
\epsilon_{t_{n,i}}^2 -\frac{1}{K}\sum_{i=n-K+1}^{n}
\epsilon_{t_{n,i}}^2 $.
We now can see the impact of Condition \ref{cond:ace2}.
To wit,
from equation (\ref{temp1a}),
\begin{eqnarray}
\XXhat^{(n)} &=&
\sum_{i=1}^{M} a_i [X, X]^{(n,K_i)}+ 
\underbrace{2\sum_{i=1}^{M} \frac{a_i}{K_i} \sum_{j=0}^{n} \epsilon_{t_{n,j}}^2 }_{=0}
+ \sum_{i=1}^{M} a_i U_{n,K_i} + \sum_{i=1}^{M} a_i V_{n,K_i} -  2 E\epsilon^2
+O_p(n^{-1/2})\notag \\
&=&\sum_{i=1}^{M} a_i [X, X]^{(n,K_i)} + \sum_{i=1}^{M} a_i U_{n,K_i} + R_n
+O_p(n^{-1/2}) ,
\label{eq:approx2}
\end{eqnarray}
where $R_n$ is the overall
remainder term, $R_n= \sum_{i=1}^{M} a_i V_{n,K_i} -  2 E\epsilon^2$.
Thus, apart from the contribution of this remainder term,
Condition \ref{cond:ace2} removes the
bias term due to $\sum \epsilon_{n,i}^2$, not only in expectation, but
almost surely. We emphasize this to stress that though we have
assumed that the $\epsilon_{t_{n,i}}$ are i.i.d., our estimator is
quite robust to the nature of the noise. 
As before,
Condition 1 assures that the first term in (\ref{eq:approx2}) will
be asymptotically unbiased for $\XX$.

Furthermore, for $i\ne l$, the $U_{n,K_i}$ and $U_{n,K_l}$ are uncorrelated. Since
  $U_{n,K_i}$ and $U_{n,K_l}$  are also the end points of zero-mean martingales, they are 
asymptotically independent as $n \rightarrow \infty$. Finally, the last term
$R_n$ is treated
separately in the proof of Theorem \ref{thm:MSRV-clt}. For now, we focus on the terms other than the $V_{n, K_i}$'s. 

If one presupposes Condition \ref{cond:ace2},  and that $R_n$ is
comparatively small, it is as if we
observe 
$$[X, X]^{(K_i)} + U_{n,K_i}, ~ i=1, ..., M.$$
Under the ideal
world of continuous observations (that is, if we
take $[X, X]^{(K_i)}$ to stand in for $\XX$), Condition \ref{cond:ace2}
makes it possible that we get $M$ (almost) independent measurements of $\XX$. This motivates the form of the MSRV estimator.

Our aim is to use Conditions \ref{cond:ace1}-\ref{cond:ace2} to construct optimal
weights $a_i$. We proceed to investigate what happens if
we just take $[X, X]^{(K_i)} \approx \XX$ in Section 3-4. From Section 5
on, we consider the more exact calculation that follows from $[X, X]^{(K_i)}  = \XX
+ O_p( (n/K_i)^{-1/2})$.

\section{Asymptotics for the Noise Term}

As above, to get a meaningful asymptotics, we let all 
quantities depend
on $n$, thus $a_i = a_{n,i}$, $M= M_n$, $K_i=K_{n,i}$, $[Y, Y]^{(K)} = [Y, Y]^{(n,K)} $, etc. 
Sometimes the dependence on $n$ is suppressed in the notation. All results are proved in Section \ref{sec:proofs}.

Consider first the noise term
\begin{equation}
\label{eq:noise}
\zeta_n = \sum_{i=1}^{M_n} a_{n,i} U_{n,K_{n,i}} 
\end{equation}
The variance of $\zeta_n$ is as follows.

\begin{proposition}
\label{prop:noisevar}
(Variance of the noise term.). Set 
$\gamma^2_n = 4 \sum_{i=1}^{M_n} {(\frac{a_{n,i}}{K_{n,i}})}^2 $. Suppose that
the $\epsilon_{t_{n,i}}$ are iid, with mean zero and $E\epsilon^2 < \infty$, and that $M_n = o(n)$ as $n \rightarrow \infty$. Then
\begin{equation}
\label{eq:varzeta}
Var(\zeta_n) = \gamma^2_n
n{(E\epsilon^2)}^2 (1 + o(1)).
\end{equation}
Also, $\gamma_n^2$ is minimized, subject to Conditions
\ref{cond:ace1}-\ref{cond:ace2}, by choosing
\begin{equation}
a_{n,i} = \frac{K_{n,i}(K_{n,i} - \bar{K_n})}{M_n Var(K_n)}
\end{equation}
where $\bar K_n=  \frac{1}{M_n} \sum_{i=1}^{M_n} K_{n,i}$ and $Var(K_n) =
\frac{1}{M_n} \sum_{i=1}^{M_n} K_{n,i}^2 - {( \frac{1}{M_n} \sum_{i=1}^{M_n}
K_{n,i})}^2$. The
resulting minimal value of $\gamma_n$ is
\begin{equation}
{\gamma^*}^{2}_n = \frac{4}{{M_n} Var(K_n)}.
\end{equation}
\end{proposition}
Since the $U_{n,K}$ are end points of martingales, by the martingale central limit theorem (\cite{hallheyde}, Chapter 3), we obtain more precisely
the following:
\begin{theorem}
\label{th:noise}
Suppose that
the $\epsilon_{t_{n,i}}$ are iid, with
$E\epsilon^2 < \infty$, and that $M = M_n = o(n)$
as $n \rightarrow \infty$. Suppose that
$\max_{1 \le i \le M_n} | a_{n,i}/(i\gamma_n) | \rightarrow 0$
as $n \rightarrow \infty$.
Then $\zeta_n/(n^{1/2}\gamma_n) \rightarrow
N(0,E(\epsilon^2)^2)$ in law, both unconditionally and conditionally on 
${\cal X}$.
\End
\end{theorem}

\section{A Class of Estimators, and Further Asymptotics
for the Noise Term \label{sec:class}}

We here develop a class of weights $a_{n,i}$ which we shall 
use in the rest of the paper. The precise form of the weights
is given in Theorem \ref{th:hft}. The rest of this section
is motivation.

In the following and for the rest of the paper, 
assume that all scales
$i=1,...,M$ are used, which is to say that $K_{n,i}=i$.
In this case, $\bar K_n = (M_n+1)/2$ and $Var(K_n)=(M_n^2 -1)/12$, and
the optimal weights from Proposition \ref{prop:noisevar} are then given by
\begin{eqnarray}
\label{eq:approxweight}
a_{n,i} &=& 12 \frac{i}{M_n^2} 
\frac{ \left ( \frac{i}{M_n} - \frac{1}{2} - \frac{1}{2M_n} \right ) }
{ \left ( 1 - \frac{1}{M^2_n} \right ) }
\end{eqnarray}
The minimum variance is given through ${\gamma^*}^{2}_n=  {48}/{[M_n(M_n^2-1)]} $,
so that 
$$Var(\zeta_n) =  {48n}{(E\epsilon^2)}^2/{[M_n(M_n^2-1)]}.$$
The form (\ref{eq:approxweight})
motivates us to consider weights on the form
\begin{equation}
\label{eq:wfunction}
a_{n,i} = \frac{1}{M_n} w_{M_n} (\frac{i}{M_n} ), ~ i=1, ..., M_n ,
\end{equation}
as this gives rise to a tractable class of estimators. We
specifically take:
\begin{equation}
\label{eq:wfunction2}
w_M(x) = xh(x) + M^{-1} x h_1(x) + 
M^{-2} x h_2(x) + M^{-3} x h_3(x)
+ o(M^{-3}) ,
\end{equation}
where $h$ and $h_1$ are functions independent of $M$. The reason for
considering this particular functional form, where $w_M(x)$ must suitably
vanish at zero, is that condition (\ref{cond:ace2}) translates 
roughly into a requirement that $\int_0^1 \frac{w_M(x)}{x} dx$ be
approximately zero.

In terms of conditions on the function $h$, Conditions
(\ref{cond:ace1})-(\ref{cond:ace2}) imply that we have to make
the following requirements on $h$:
\begin{condition}\label{cond:hearts1}
$\int_0^1 x h(x) dx =1,$
\end{condition}
\begin{condition}\label{cond:hearts2}
$\int_0^1 h(x) dx = 0.$
\end{condition}
With  slightly stronger requirement on $h$, we can show that (\ref{eq:approx2}) holds more generally.

\begin{theorem}
\label{th:hft}
Let $h_0 = h$, and suppose that for $i = 0, ..., 2$,
$h_i$ is $3-i$ times continuously differentiable on $[0,1]$,
and that $h_3$ is continuous on $[0,1]$. Suppose
that $h$
satisfies Conditions \ref{cond:hearts1}-\ref{cond:hearts2}.
Also assume that
\begin{eqnarray}
\int_0^1 h_1(x)dx
+ \frac{1}{2} ( h(1) - h(0) ) &=& 0 , \notag \\
\int_0^1 h_2(x)dx
+ \frac{1}{2} ( h_1(1) - h_1(0)  )
+ \frac{1}{12} ( h'(1) - h'(0) )  &=& 0 ,
\label{eq:extracond}
\\
\mbox{and } \int_0^1 h_3(x)dx
+ \frac{1}{12} ( h_1'(1) - h_1'(0)  ) &=& 0 .
\notag 
\end{eqnarray}
Let the
$a_{n,i}$ be given by (\ref{eq:wfunction})-(\ref{eq:wfunction2}),
where the $o(M^{-3})$ is uniform in $x \in [0,1]$.
Finally, suppose that the $\epsilon_{t_{n,i}}$ are i.i.d., with
$E\epsilon^2 < \infty$. 
Then approximation (\ref{eq:approx2}) remains valid, up
to $o_p(n/M^3_n)$.
\end{theorem}

{\bf The final class of estimators.}
Our estimation procedure will in the following be using weights $a_{n,i}$
which satisfy the description in Theorem \ref{th:hft}.
\begin{remark}\label{remark.thm2}
[Comments on Theorem \ref{th:hft}:] By adding terms in
(\ref{eq:wfunction2}), one can make the approximation 
in (\ref{eq:approx2}) as good as one wants (up to $O_p(n^{-1/2})$).
We will later use $M_n = O(n^{1/2})$, which is why we have
chosen the given number of terms in (\ref{eq:wfunction2}). Also, it should be noted that
the approximation to Condition \ref{cond:ace2} has to be much
finer than to Condition \ref{cond:ace1}, since we
are seeking to make $\sum_{i=1}^{M} \frac{a_i}{K_i} \sum_{i=0}^{n} \epsilon_{t_{n,i}}^2=
n \left( \sum_{i=1}^{M} \frac{a_i}{K_i} \right ) E\epsilon^2 (1 + o_p(1))$
negligible for asymptotic purposes.
\end{remark}
As we shall see, the specific choices for 
$h_1$, $h_2$, and $h_3$ do not play any role in any of the
later expressions for asymptotic variance.
A simple choice of $h_1$ which satisfies (\ref{eq:extracond})
is given by $h_1(x) = - h^\prime(x)/2$, with $h_2(x) = h_2$ and
$h_3(x) = h_3$, both constants. In this case, $h_2 = - (h'(1)-h'(0))/6$
and $h_3 = (h''(1)-h''(0))/24$. With this choice, one obtains
\begin{equation}
\label{eq:hfunction}
a_{n,i} = \frac{i}{M_n^2} h(\frac{i}{M_n} ) - \frac{1}{2}
\frac{i}{M_n^3} h^\prime (\frac{i}{M_n} )
+ \frac{i}{M^3_n} h_2 + \frac{i}{M^4_n} h_3
\end{equation}
For the noise-optimal weights in
(\ref{eq:approxweight}) at the end of Section 3, 
$h$ takes the form
\begin{equation}
\label{eq:hstar}
h^*_\zeta (x) = 12 \left (x - \frac{1}{2} \right ) .
\end{equation}
Under this choice, the $a_{n,i}$ given by (\ref{eq:hfunction}) is identical 
to the one in (\ref{eq:approxweight}), up to a negligible multiplicative factor of 
$( 1 - M^{-2}_n)^{-1}$.

The following corollary to Theorem \ref{th:noise} is now
immediate, since
$\gamma^2_n = 4 M_n^{-3} \int_0^1 {h(x)^2} dx (1 + o(1))$
as $n \rightarrow \infty$.
\begin{corollary}
\label{cor:noise}
Suppose that
the $\epsilon_{t_{n,i}}$ are iid, with
$E\epsilon^2 < \infty$, and that $M = M_n = o(n)$
as $n \rightarrow \infty$. Also assume that
the $a_{n,i}$ are given by (\ref{eq:wfunction}), and that
the conditions of Theorem \ref{th:hft} are satisfied.
Then $(M_n^{3}/n)^{1/2}\zeta_n \rightarrow
N(0,4 E(\epsilon^2)^2 \int_0^1 {h(x)^2} dx $ in law, both unconditionally 
and conditionally on ${\cal X}$.
\End
\end{corollary}

\section{Asymptotics of the Discretization Error \label{sec:as-dis}}

We have obtained the optimal weights as far as reducing the noise
is concerned. However, as in (\ref{eq:approx2}), there remains
two types of error: the {\it discretization error},
due to the fact that the observations only take
place at discrete time points, along with the {\it residual} $R_{n}$,
which also will turn out to not quite vanish. We study these in turn,
and then state a result for the total asymptotics for the MSRV estimator.

For the discretization error, we need some additional concepts. 
\begin{definition}
Let $0 = t_{n,0} < t_{n,1} < ... < t_{n,n} = T$ be the observation
times when there are $n$ observations. We refer to
$\G_n = \{t_{n,0} , t_{n,1} , ... , t_{n,n}\}$ as a ``grid''
or a ``partition''
of $[0,T]$.
Following Section 2.6 of \cite{myklandzhang}, the ``Asymptotic Quadratic 
Variation of Time'' ("AQVT") $H(t)$ is defined by
\begin{equation}
\label{eq:qvtime}
H(t) = \lim_{n \rightarrow \infty} \frac{n}{T} \sum_{t_{n,i+1} \le t}
(t_{n,i} - t_{n,i-1})^2 ,
\end{equation}
provided the limit exists.
\end{definition}
We  assume that
\begin{equation}
\max_{1 \le i \le n}~|t_{n,i+1}-{t_{n,i}}|~~=O\left(  \frac{1}{n}\right)  , 
\label{eq:maxdelta}
\end{equation}
whence every subsequence has a subsequence so that the asymptotic
quadratic variation of time exists. From an applied point of view,
 there is little loss in assuming the existence of
the asymptotic quadratic variation of time, cf.
the argument
at the very end of \cite{twoscales} (on p. 1411).

Note that from (\ref{eq:maxdelta}), $H(t)$ is Lipschitz continuous
provided it exists. We give the following change-of-variable rule
for the AQVT:

\begin{lemma}
\label{lemma:ch-var}
(Change of variables in the AQVT.)
Assume (\ref{eq:maxdelta}) and that the AQVT H(t) exists. Let
$G: [0,T] \rightarrow [0,T]$ be Lipschitz continuous. Set 
$u_{n,i} = G(t_{n,i})$. Then
\begin{equation*}
K(u) = \lim_{n \rightarrow \infty} \frac{n}{T} \sum_{u_{n,i} \le u}
(u_{n,i} - u_{n,i-1})^2 
\end{equation*}
exists, and
\begin{equation}
\label{eq:ch-var}
H'(t) G'(t) = K'(G(t))
\end{equation}
almost everywhere on $[0,T]$.
\end{lemma}

The following result is also useful and illustrative.
\begin{lemma}
\label{lemma:approxpoints}
Assume the conditions of Lemma \ref{lemma:ch-var}. Then
$K(T) = T$ if and only if
\begin{equation}
\label{eq:approxpoints}
\sum_{i=0}^{n} \left ( u_{n,i} - u_{n,i-1} - \frac{T}{n} \right )^2 = o(n^{-1}) .
\end{equation}
\end{lemma}

\begin{remark}
\label{remark:almostunif}
The importance of these two lemmas is that one can compare irregular
and ``almost equidistant'' sampling. If $H'(t)$ exists, is continuous,
and is bounded
below by a constant $c > 0$, one can define $G(t) = \int_0^t H'(s)^{-1}ds$,
and consider the process $\tilde X_u = X_{G(u)}$. This process
satisfies the same regularity conditions as those that we impose on
$X$, and, furthermore, the sampling times $u_{n,i} = G(t_{n,i})$ 
are close to equidistant in the sense of equation (\ref{eq:approxpoints}).
The further implication of this is discussed in Remark \ref{remark:almostunif2}
after Theorem \ref{thm:discrerror}.
\end{remark}

Define $\eta$ as the nonnegative square root of
\begin{equation}
\eta^2 = \int_0^T H^\prime (t) \sigma_t^4 dt
\end{equation}
Finally, we define ``stable convergence''. 
\begin{definition} 
If $Z_n$ is a sequence of ${\cal X}$-measurable random variables,
the $Z_n$ converges stably in law to $Z$ as $n \rightarrow \infty$ 
if there is an extension
of ${\cal X}$ so that for all $A \in {\cal X}$ and for all
bounded continuous $g$, $EI_A g(Z_n) \rightarrow EI_A g(Z)$ as
$n \rightarrow \infty$.
\end{definition}
For further discussion of
stable convergence, see \cite{renyi63}, \cite{aldouseagleson78},
Chapter 3 (p. 56) of \cite{hallheyde}, \cite{rootzen80} and Section 2
(p. 169-170) of \cite{jacodprotter98}.
It is a useful device in operationalizing asymptotic
conditionality. There is some choice in what one takes as the $\sigma$-field
${\cal X}$ in this definition.

We can now state the main theorem for the asymptotic
behavior of finitely many of the $[X, X]^{(K)}=[X, X]^{(n,K)}$.
\begin{theorem}
\label{thm:discrerror}
(CLT for the discretization error in $[X, X]^{(K)}$.) Suppose the structure of $X$ follows Assumption \ref{ass:process}. Also suppose that
the observation times $t_{n,i}$ are nonrandom, satisfy 
(\ref{eq:maxdelta}), and that the asymptotic quadratic variation of time
$H(t)$ exists and is continuously differentiable. 
Assume that $\min_{0 \le t \le T} H'(t) > 0$. 
Let $M_n \rightarrow \infty$ as $n \rightarrow \infty$, with
$M_n = o(n)$.
Let $(K_{n,1}, ..., K_{n,L}) / M_n \rightarrow
(\kappa_1 , ..., \kappa_L)$ as $n \rightarrow \infty$. Let
$\Gamma$ be an $L \times L$ matrix
with $(I,J)$ entry given by
\begin{equation}
\label{eq:Gamma}
\Gamma_{I,J} = \frac{2}{3} T \min(\kappa_I,\kappa_J) 
\left (3 - \frac{\min(\kappa_I,\kappa_J)}{\max(\kappa_I,\kappa_J)} \right ),
\end{equation}
and let $Z$ be a
normal random vector with covariance matrix $\Gamma$.
Let $Z$ be independent of ${\cal X}$. Then, as
$n \rightarrow \infty$
the vector 
$(n/M_n)^{1/2}([X, X]^{(n,K_{n,1})} - \XX , ..., [X, X]^{(n,K_{n,L})} - \XX )$
converges stably in law to $\eta Z$.
\end{theorem}

\begin{remark}
\label{remark:almostunif2}
Even in the scalar ($L=1$) case, this result in Theorem
\ref{thm:discrerror} is a gain over our earlier Theorem
3 (p. 1401) in \cite{twoscales}. To characterize the asymptotic 
distribution we use an asymptotic quadratic variation of time (AQVT)
which is independent of choice of scale and coincides with the
original object introduced in \cite{myklandzhang} (Section 2.6).
This is unlike the time variation measure used in section 3.4
in \cite{twoscales}, and Theorem \ref{thm:discrerror} provides
a substantial simplification
of the asymptotic expressions. To do this, we have used the
approach described above in Remark \ref{remark:almostunif}.

It is conjectured that the regularity conditions for Theorem \ref{thm:discrerror}  can
be reduced to those of Proposition 1 of \cite{myklandzhang},
but investigating this is beyond the scope of this paper.
\end{remark}
As a corollary to Theorem \ref{thm:discrerror}, we now finally obtain the
asymptotics for the discretization part of the MSRV, as follows.
\begin{corollary}
\label{cor:discrerror}
(CLT for the discretization error in the MSRV.) Let $a_{n,i}$ satisfy (\ref{eq:wfunction})-(\ref{eq:wfunction2}), and let the
conditions of Theorem \ref{th:hft} be satisfied.
Further, make Assumption \ref{ass:process}. Also suppose that
the observation times $t_{n,i}$ are nonrandom, satisfy
(\ref{eq:maxdelta}), and that the asymptotic quadratic variation of time
$H(t)$ exists and is continuously differentiable.
Assume that $\min_{0 \le t \le T} H'(t) > 0$.
Let $M_n \rightarrow \infty$ as $n \rightarrow \infty$, with
$M_n/n = o(1)$ and $M_n^3/n \rightarrow \infty$. Set
\begin{equation}
\label{eq:discrerror1}
\eta^2_h =\frac{4}{3} T \eta^2
\int_0^1 dx \int_0^x h(y)h(x) y^2 \left( 3x - {y} \right ) dy
\end{equation}
Then
\begin{equation}
\label{eq:discrerror2}
(n/M_n)^{1/2} \left (
\sum_{i=1}^{M_n} a_{n,i} [X, X]^{(n,i)} - \XX \right ) \rightarrow 
\eta_h Z
\end{equation}
stably in law, where $Z$ is standard normal and independent of ${\cal X}$.
\end{corollary}
\begin{remark} Note that the condition $M_n^3/n \rightarrow \infty$
is present because we have not imposed too many conditions on
$h$; if it were necessary, the assumption could be removed
by considering a slightly smaller class of $h$s.
\end{remark}

\section{Overall Asymptotics
for the MSRV Estimator}

There are two main sources of error in the MSRV. 
On the one hand,
we have seen in Corollary \ref{cor:noise}
(at the end of Section \ref{sec:class}) that
if $M_n$ time scales are used, the part of $\XXhat^{(n)} - \XX$
which is due purely to the noise $\epsilon$
can be reduced to have order $O_p(n^{1/2} M_n^{-3/2})$. At the same
time, Corollary \ref{cor:discrerror} shows that
the pure discretization error is of order $O_p(n^{-1/2} M_n^{1/2})$.
To balance these two terms, the optimal $M_n$ is therefore of the order
\begin{equation}
M_n = O(n^{1/2}) ,
\end{equation}
assuming that the remainder term in (\ref{eq:approx2}) does not
cause problems, which is indeed the case.
This leads to a variance-variance tradeoff, and
the rate of convergence for the MSRV estimator is then
$\XXhat^{(n)} - \XX = O_p(n^{-1/4})$.
This result is an improvement on the two scales estimator, for which
the corresponding rate is $O_p(n^{-1/6})$. We embody this in
the following result.

\begin{theorem}
\label{thm:MSRV-clt}
Let $a_{n,i}$ satisfy (\ref{eq:wfunction})-(\ref{eq:wfunction2}), and let the
conditions of Theorem \ref{th:hft} be satisfied.
Further, make Assumptions \ref{ass:process}-\ref{ass:noise}. Also suppose that
the observation times $t_{n,i}$ are nonrandom, satisfy
(\ref{eq:maxdelta}), and that the asymptotic quadratic variation of time
$H(t)$ exists and is continuously differentiable.
Assume that $\min_{0 \le t \le T} H'(t) > 0$. Suppose that
$M_n/n^{1/2} \rightarrow c$ as $n \rightarrow \infty$.
Let $Z$ be a standard normal random variable independent of
${\cal X}$. Set
\begin{eqnarray}
\nu^2_h&=&4 c^{-3} {(E\epsilon^2)}^2 \int_0^1 h(x)^2 dx
+ c \frac{4}{3} T \eta^2
\int_0^1 dx \int_0^x h(y)h(x) y^2 \left( 3x - y \right ) dy\notag\\
& +&\hspace*{-0.1in}4 c^{-1} Var(\epsilon^2) \int_0^1 \int_0^y
x h(x)h(y) dx dy  + 8 c^{-1} E\epsilon^2 \int_0^1 \int_0^1 h(x)h(y) {\min(x, y)} dx dy \XX ~~ .
\end{eqnarray}
Then
\begin{equation}
n^{1/4} \left ( \XXhat^{(n)} - \XX \right ) \rightarrow \nu_h Z ,
\end{equation}
stably in law, as $n \rightarrow \infty$.
\end{theorem}

For the noise optimal $h$-function from equation (\ref{eq:hstar})
(cf. equation  (\ref{eq:approxweight})), we can now
calculate 
the value of the asymptotic variance of the MSRV. Note that
if $h(x) = 12(x-1/2)$, we obtain
\begin{eqnarray*}
\int_0^1 dx \int_0^x h(y)h(x) y^2 \left( 3x - y \right ) dy =
\frac{39}{35},\\
\int_0^1 \int_0^y
x h(x)h(y) dx dy =\frac{3}{5},\\
\int_0^1 \int_0^1 h(x)h(y)
{\min(x, y)} dx dy =\frac{6}{5}.
\end{eqnarray*}
Hence, in this case, the asymptotic variance becomes
\begin{equation}
\nu^2_h=48 c^{-3} {(E\epsilon^2)}^2
+ \frac{52}{35} c T \eta^2  + \frac{12}{5} c^{-1} Var(\epsilon^2)
 +  \frac{48}{5} c^{-1} E\epsilon^2 \XX
\end{equation}

\section{Conclusion}
In this paper, we have introduced the {\it Multi Scale Realized
Volatility (MSRV)} and shown a central limit theorem (Theorem
\ref{thm:MSRV-clt}) for this
estimator. This permits the setting of intervals for the true
integrated volatility on the basis of the MSRV. As a consequence
of our result, it is clear that the MSRV is rate efficient, with
a rate of convergence of $O_p(n^{-1/4})$.

In terms of the general study of realized volatilities, Section
\ref{sec:as-dis} also shows further properties of the asymptotic
quadratic variation of time (AQVT), as earlier introduced by
\cite{myklandzhang} and \cite{twoscales}. In particular, Theorem
\ref{thm:discrerror} shows that one can use the regular one-step
AQVT also for multistep realized volatilities, thus improving
on Theorems 2 and 3 (p. 1401) in \cite{twoscales}.

Finally, note that most of the arguments we have used hold up
also when the noise process $\epsilon_{t_{n,i}}$ is no longer iid.
One can, for example, model this process as being stationary
(but with mean zero). If the process is sufficiently mixing,
this will change the asymptotic variance of the MSRV, but not the consistency, nor 
the convergence rate of $O_p(n^{-1/4})$, see for example Chapter
5 of \cite{hallheyde} for the basic limit theory for dependent
sums. However, we have not sought to develop the specific conditions
for the CLT to hold in the case when the process is mixing.

\section{Proofs of Results. \label{sec:proofs}}

Note that for ease of notation, we sometimes suppress the dependence on $n$  in the notation. For example, $a_i = a_{n,i}$, $M= M_n$, $K_i=K_{n,i}$, $[Y, Y]^{(K)} = [Y, Y]^{(n,K)} $, etc. 
Also, we in this section write $t_i$ for $t_{n,i}$, to avoid cluttering of the notations. 

\subsection{Proof of Proposition \ref{prop:approxA}.}
Write
\begin{eqnarray}
\XXhat^{(n)} &=& \sum_{i=1}^{M} a_i [Y, Y]^{(n,K_i)} +
(\alpha_1 - a_1) ( [Y, Y]^{(n,K_1)} - [Y, Y]^{(n,K_2)}) \notag \\
&=& \sum_{i=1}^{M} a_i [Y, Y]^{(n,K_i)} - 2 E \epsilon^2 +
O_p(n^{-1/2})\label{temp1}
\end{eqnarray}
where the final approximation follows from Lemma 1 (p. 1398)
in \cite{twoscales}.

\subsection{Proof of Proposition \ref{prop:noisevar}}

Since $U_{n,K_{n,i}}$ and $U_{n,K_{n,l}}$ are uncorrelated
($i \ne l$) zero-mean 
martingales,
\begin{eqnarray}
Var(\zeta_n )&=&
\sum_{i=1}^{M_n} a_{n,i}^2 Var(U_{n,K_{n,i}} ) \notag \\
&=&4 \sum_{i=1}^{M_n} {(\frac{a_{n,i}}{K_{n,i}})}^2
(n-K_{n,i} +1){(E\epsilon^2)}^2 \notag \\
&=& \gamma^2
n{(E\epsilon^2)}^2 (1 + o(1)),
\label{eq:varzeta2}
\end{eqnarray}
showing equation (\ref{eq:varzeta}).
The last transition in (\ref{eq:varzeta2}) follows because
$M_n = o(n)$.

We minimize $\gamma_n^2$, subject to the constraints in Conditions 
\ref{cond:ace1}-\ref{cond:ace2}. This is
established by setting
\begin{eqnarray*}
\frac{\partial}{\partial a_{n,i}} [\gamma^2_n + \lambda_1 (\sum a_{n,i} -1) + \lambda_2 (\sum
  \frac{a_{n,i}}{K_{n.i}} )]=8 \frac{a_{n,i}}{K_{n,i}^2} +\lambda_1 +
  \frac{\lambda_2}{K_{n,i}}
\end{eqnarray*}
to zero, resulting in  $a_{n,i} = -\frac{1}{8} (\lambda_1 K_{n,i}^2+ \lambda_2 K_{n,i})$. 
One can determine the $\lambda$'s by solving
$$
\left \{
\begin{array}{l}
1=\sum_{i=1}^{M_n} a_{n,i} = -\frac{1}{8} (\lambda_1 \sum_{i=1}^{M_n} K_{n,i}^2+
  \lambda_2 \sum_{i=1}^{M_n} K_{n,i})\\
0=\sum_{i=1}^{M_n} \frac{a_i}{K_{n,i}}= -\frac{1}{8} (\lambda_1 \sum_{i=1}^{M_n}
  K_{n,i}+ M_n \lambda_2 ) 
\end{array}
\right.
$$
This leads to 
$$\lambda_1 = -\frac{8}{M_n Var(K_n)} \quad \mbox{ and } \quad \lambda_2= 
\frac{8 \bar K_n}{M_n  Var(K_n)},$$
where $\bar K_n$ and $Var(K_n)$ are as given in Proposition
(\ref{prop:noisevar}). This shows the rest of the proposition.

\subsection{Proof of Theorem \ref{th:noise}.}
Assume without loss of generality that $K_i=i$ for $i=1, ..., M$.
To avoid cluttering the notation, we write $a_i$ for $a_{n,i}$.
Note that $\zeta_n$ is the end point of a martingale. We show that
$\zeta_n/(n^{1/2}\gamma_n)$ satisfies the conditions of the version
of the Martingale Central Limit Theorem which is stated in
Corollary 3.1 (p. 58-59) of \cite{hallheyde}. The result then
follows. Note that we shall take, in the notation of
\cite{hallheyde}, ${\cal F}_{n, j}$ to be the smallest $\sigma$-field
making $\epsilon_{t_i}$, $i= 1, ... , j$, and the whole $X_t$ process,
measurable.

We start with the Lindeberg condition. For given $\delta$, define
$f_\delta(x) = E(\epsilon^2 x^2 I_{\{ | \epsilon x | > \delta \} })$.
Also set 
\begin{equation*}
r_n(x) = Ef_{\delta n^{1/2}} \left (
- \frac{1}{\gamma_n} \sum_{i=1}^{M_n \wedge j} \frac{2 a_i}{i} \epsilon_{t_i} \right )
\mbox{   for } \frac{j-1}{n} \le x < \frac{j}{n} .
\end{equation*}
We then obtain
\begin{align}
\label{eq:Lindeberg1}
\sum_{j=1}^n E &\left ( \epsilon_{t_j}^2 
\left (
- \frac{1}{n^{1/2} \gamma_n} \sum_{i=1}^{M_n \wedge j} \frac{2 a_i}{i} \epsilon_{t_{j-i}} \right )^2
I_{\{ | \epsilon_{t_j}
\left (
- \frac{1}{n^{1/2} \gamma_n} \sum_{i=1}^{M_n \wedge j} \frac{2 a_i}{i} \epsilon_{t_{j-i}} \right ) |
> \delta \} } \right ) \notag \\
&= \frac{1}{n} \sum_{j=1}^n
E f_{\delta n^{1/2}} \left (
- \frac{1}{\gamma_n} \sum_{i=1}^{M_n \wedge j} \frac{2 a_i}{i} \epsilon_{t_{j-i}} \right ) \notag \\
& = \int_0^1 r_n(x) dx \mbox{ since the }\epsilon_{t_i}\mbox{ are i.i.d.} \notag \\
& \rightarrow 0 \mbox{ as } n \rightarrow \infty,
\end{align}
where the last transition is explained in the next paragraph. By
Chebychev's inequality, the conditional Lindeberg condition in
Corollary 3.1 of \cite{hallheyde} is thus satisfied.

The last transition in (\ref{eq:Lindeberg1}) is because of the following. 
First fix $x \in [0,1)$,
and let $j_n$ be the corresponding $j$ in the definition of $r_n(x)$. 
Let $Z_n = - \frac{1}{\gamma_n} \sum_{i=1}^{M_n \wedge j_n} \frac{2 a_i}{i} \epsilon_{t_i} $, so that $r_n(x) = Ef_{\delta n^{1/2}} (Z_n)$.

Note that $Z_n$ is a sum of independent random variables which,
satisfies the Lindeberg condition:
\begin{equation*}
\sum_{i=1}^{M_n \wedge j_n} E\left ( \frac{- 2 a_i}{i\gamma_n} \epsilon_{t_i}
\right )^2 I_{\{ | \frac{- 2 a_i}{i\gamma_n} \epsilon_{t_i} | > \delta \} }
= \sum_{i=1}^{M_n \wedge j_n} f_\delta \left ( \frac{- 2 a_i}{i\gamma_n} \right ) \rightarrow 0
\end{equation*}
as $n \rightarrow \infty$, since $\max_{i} | a_i/i\gamma_n | \rightarrow 0$.
The ensuing asymptotic normality of $Z_n$ (if necessary by going to 
subsequences of subsequences) shows that $r_n(x) \rightarrow 0$
as $n \rightarrow \infty$. Since $0 \le r_n(x) \le 1$, the final
transition in (\ref{eq:Lindeberg1}) follows by dominated convergence.

We now turn to the sum of conditional variances in the corollary in
\cite{hallheyde}.
\begin{align}
\sum_{j=1}^n E & \left ( \epsilon_{t_j}^2
\left (
- \frac{1}{n^{1/2} \gamma_n} \sum_{i=1}^{M_n \wedge j} \frac{2 a_i}{i} \epsilon_{t_{j-i}} \right )^2 | {\cal F}_{n, j-1} \right ) \notag \\
& = E (\epsilon^2) \frac{1}{n \gamma_n^2}
\sum_{j=1}^n \left ( \sum_{i=1}^{M_n \wedge j} 
\frac{2 a_i}{i} \epsilon_{t_{j-i}} \right )^2 \notag \\
& = 1 + o_p(1) .
\label{eq:sumcondvar}
\end{align}
The last transition is obvious by appealing to M-dependence.
A rigorous but tedious proof is obtained by splitting the sum
into main terms of the type $\epsilon^2_{t_i}$ and cross-terms of the form
$\epsilon_{t_i} \epsilon_{t_j}$ ($i \ne j$). 

In view of (\ref{eq:Lindeberg1})-(\ref{eq:sumcondvar}), Theorem
 \ref{th:noise} is proved by using Corollary 3.1 and the Remarks
following this corollary (p. 58-59) in \cite{hallheyde}.

\subsection{Proof of Theorem \ref{th:hft}.} 
We need to show that 
$\sum_{i=1}^{M} \frac{a_{n,i}}{K_{n,i}} \sum_{i=0}^{n} \epsilon_{t_i}^2=o_p(n/M^3_n)$,
in other words, we need $\sum_{i=1}^{M-n} \frac{a_{n,i}}{K_{n,i}} = o(M_n^{-3})$.
By Taylor expansion
\begin{eqnarray}
\frac{1}{M} \sum_{i=1}^{M} h(\frac{i}{M} )
&=& \int_0^1 h(x) dx + \frac{1}{2 M^2} \sum_{i=1}^{M} h'(\frac{i}{M} )
- \frac{1}{3! M^3} \sum_{i=1}^{M} h''(\frac{i}{M} )
+ \frac{1}{4! M^4} \sum_{i=1}^{M} h'''(\frac{i}{M} ) +o(M^{-3})
\notag \\
&=& \int_0^1 h(x) dx + \frac{1}{2 M} (h(1)-h(0))
+ \frac{1}{12 M^3} \sum_{i=1}^{M} h''(\frac{i}{M} )
- \frac{1}{24 M^4} \sum_{i=1}^{M} h'''(\frac{i}{M} ) +o(M^{-3})
\notag \\
&=& \int_0^1 h(x) dx + \frac{1}{2 M} (h(1)-h(0))
+ \frac{1}{12 M^2} (h'(1)-h'(0)) +o(M^{-3}) ,
\end{eqnarray}
where the later line follows by iterating the first line.
By similar argument on $h_1$ to $h_3$,
\begin{eqnarray}
\frac{1}{M} \sum_{i=1}^{M} (\frac{i}{M})^{-1} w_{M} (\frac{i}{M} ) &=& 
\frac{1}{M} \sum_{i=1}^{M} h(\frac{i}{M} )
+ \frac{1}{M^2} \sum_{i=1}^{M} h_1(\frac{i}{M} ) 
+ \frac{1}{M^3} \sum_{i=1}^{M} h_2(\frac{i}{M} )
+ \frac{1}{M^4} \sum_{i=1}^{M} h_3(\frac{i}{M} ) +o(M^{-3}) \notag \\
&=& \int_0^1 h(x)dx + \frac{1}{M} \left ( \int_0^1 h_1(x)dx
+ \frac{1}{2} ( h(1) - h(0) ) \right ) \notag \\
&&+ \frac{1}{M^2} \left ( \int_0^1 h_2(x)dx
+ \frac{1}{2} ( h_1(1) - h_1(0)  )
+ \frac{1}{12} ( h'(1) - h'(0) ) \right ) \notag \\
&&+ \frac{1}{M^3} \left ( \int_0^1 h_3(x)dx
+ \frac{1}{12} ( h_1'(1) - h_1'(0)  ) \right )
+ o(\frac{1}{M^3}) \notag \\
&=& o(\frac{1}{M^3}), \notag
\end{eqnarray}
by (\ref{eq:extracond}). This shows the result.

\subsection{Proof of Lemma \ref{lemma:ch-var}.}
To get the rigorous statement, we proceed as follows. 
Every subsequence has a
further subsequence for which $K(u)$ exists, and this $K$
is obviously Lipschitz continuous. We will show that
(\ref{eq:ch-var}) hold. Since this equation is independent of subsequence,
the result will have been proved.

Let $B_t$ be a standard Brownian motion, 
and let 
$\tilde{B}_t = B_{G(t)}$. By comparing the asymptotic
distributions of $(T/n)^{-1/2} [ \sum_{t_{i} \le t }
(\tilde B_{t_{i}} - \tilde B_{t_{i-1}})^2 - <\tilde B, \tilde B>_t ]$
and 
$(T/n)^{-1/2} [ \sum_{u_{i} \le u }
(B_{u_{i}} - B_{u_{i-1}})^2 - <B,B>_u ]$, we obtain from Proposition 1
of \cite{myklandzhang} that
\begin{equation*}
\int_0^t 2 H'(s) (<\tilde B, \tilde B>_s')^2 ds
= \int_0^{G(t)} 2 K'(v) (<B,B>_v')^2 dv \quad \mbox{ for all }
t \in [0.T] .
\end{equation*}
Since $<B,B>_v' = 1$ and $<\tilde B, \tilde B>_s' = G'(s)$ a.e., 
equation (\ref{eq:ch-var}), and hence the lemma, follows.

\subsection{Proof of Lemma \ref{lemma:approxpoints}.}
Set $\delta_{n,i} = u_{n,i} - u_{n,i-1} - T/n$.
Then
\begin{eqnarray*}
\frac{n}{T} \sum_{i}
\left ( u_{n,i} - u_{n,i-1} \right )^2 &= \frac{n}{T} \sum_{i}
\left ( \frac{T}{n} + \delta_{n,i} \right )^2 \notag \\
&= T + 2 \sum_{i} \delta_{n,i}
+ \frac{T}{n}  \sum_{i} \delta_{n,i}^2.
\end{eqnarray*}
Since $\sum_{i} \delta_{n,i} = 0$, the Lemma follows by letting 
$n \rightarrow \infty$.

\subsection{Proof of Theorem \ref{thm:discrerror}.}
Following Lemmas \ref{lemma:ch-var} and
\ref{lemma:approxpoints}, and Remark \ref{remark:almostunif},
we can assume without loss of generality that the $t_{n,i}$ satisfy
(in place of the $u_{n,i}$) the equation (\ref{eq:approxpoints}).

Consider the scalar case ($L=1$) first, with $K_{n} = K_{n,1} = M_n$.
In the sequel, all prelimiting quantities are subscripted by $n$,
and we suppress the $n$ for ease of notation (except when it
seems necessary).
We now refer to Theorems 2 and 3 (p. 1401) in \cite{twoscales}.
Use the notation $\Delta t_i$, $h_i$ and $\eta_n$ as in that paper, and let
$\overline{\Delta t} = T/n$. (Note that the usage of ``$\eta$''
in this paper is different from that of \cite{twoscales}.
Also define
\begin{equation*}
\tilde{h}_{i}=\frac{4}{K \overline{\Delta t}} \sum_{j=1}^{(K-1)\wedge i} (1-\frac
{j}{K})^{2} \overline{\Delta t}
\mbox{ and }
\tilde{\eta}_{n}^{2}~=~\sum_{i}\tilde{h}_{i}\sigma^4_{t_{i}} \overline{\Delta t} .
\end{equation*}
Note that if we show that $\tilde{\eta}_n - \eta_n \rightarrow 0$
in probability
as $n \rightarrow \infty$, we have shown the scalar version of the 
theorem. This is because we we will then have shown that the conditions
of the two Theorems in \cite{twoscales} are satisfied, and that we can
calculate the asymptotic variances as if $t_{n,i} = iT/n$.

To this end, note first that
\begin{align}
| \sum_i h_{i}\sigma^4_{t_{i}} (\Delta t_{i} - \overline{\Delta t} ) |
& \le (\sigma^+)^4 \left ( \sum_i h_{i}^2 \right )^{1/2} 
\left ( \sum_i (\Delta t_{i} - \overline{\Delta t} )^2 \right )^{1/2}
\notag \\
& = O(n^{1/2}) \times o(n^{-1/2}) = o(1) ,
\label{eq:etadiff1}
\end{align}
where the orders follow, respectively, from equation (45)
in \cite{twoscales}, and equation (\ref{eq:approxpoints}) in this paper.
Then note that
\begin{align}
| \sum_i (h_{i} - \tilde h_i) \sigma^4_{t_{i}} \overline{\Delta t} |
& = | \frac{4}{K} (\sigma^+)^4 \sum_{j=1}^{K-1} \left (1-\frac {j}{K} \right )^{2}
\left ( \sum_{l = (K-1)^+}^{n-j} ( \Delta t_l - \overline{\Delta t} )
\right ) | \notag \\
& \le \frac{4}{K} (\sigma^+)^4 \sum_{j=1}^{K-1} \left (1-\frac {j}{K} \right )^{2}
\times \left ( \sum_i (\Delta t_{i} - \overline{\Delta t} )^2 \right )^{1/2}
\notag \\
& = O(1) \times o(n^{-1/2}) = o(1) 
\label{eq:etadiff2}
\end{align}
where, again, the orders follow, respectively, from equation (45)
in \cite{twoscales}, and equation (\ref{eq:approxpoints}) in this paper.

Equations (\ref{eq:etadiff1})-(\ref{eq:etadiff2}) combine to show that
$\tilde{\eta}_n - \eta_n \rightarrow 0$
in probability
as $n \rightarrow \infty$.

For the general ($L > 1$) case, first note that 
since $\mu_t$ and $\sigma_t$ are bounded (Assumption \ref{ass:process}),
by Girsanov's Theorem (see, for example,
Chapter 3.5 (pp. 190-201) of \cite{karatzasshreve}, or
Chapter II-3b (pp. 168-170) of
\cite{jacodshiryaev2003}), we can without loss of generality
further suppose that $\mu_t=0$ identically.
This is because of the stability of the convergence, cf.
the methodology in \cite{rootzen80}.

Now set
\begin{eqnarray*}
(X,X)^{(K)} = \frac{2}{K} \sum_{j=0}^{n-1} (X_{t_{j+1}}- X_{t_j}) \sum_{r=1}^{j \wedge (K-1)}
(K-r) (X_{t_{j-r+1}}- X_{t_{j-r}})
\end{eqnarray*}
and note that 
\begin{eqnarray*}
[X, X]^{(n,K)} &=& (X,X)^{(K)} + [X, X]^{(n,1)} + O_p(K/n)  \\
&=& (X,X)^{(K)}_T + \XX + O_p(n^{-1/2})  + O_p(K/n) ,
\end{eqnarray*}
from Proposition 1 in \cite{myklandzhang}.

Let $M_t^{n,I}$ be the continuous martingale for which
$M_T^{n,I}= (X,X)^{(I)} (n/M_n)^{1/2}$. The proof of Theorem
2 in \cite{twoscales} actually
establishes that the sequence of processes $(M_t^{n,K_{n,I}})$
is C-tight in the sense of Definition VI.3.25 (p. 351) of
\cite{jacodshiryaev2003}. This is because of Theorem VI.4.13 (p. 358)
and Corollary VI.6.30 (p. 385), also in \cite{jacodshiryaev2003}.
The same corollary then establishes that asymptotic distribution
is as described in Theorem \ref{thm:discrerror}, provided we
can show that 
\begin{equation}
\label{eq:qvconv}
<M^{n,K_{n,I}},M^{n,K_{n,J}}>_T \rightarrow \eta^2 \Gamma 
\mbox{ as } n \rightarrow \infty .
\end{equation}
This is because
of L{\'e}vy's Theorem (see \cite{karatzasshreve}, Theorem 3.16, p. 157). 
The stable convergence follows as
in the proof of Theorem 3 of \cite{twoscales}, the conditions
for which have already been satisfied.

We finally need to show (\ref{eq:qvconv}). As in the scalar case,
we assume (\ref{eq:approxpoints}), and the same kind argument used in the
scalar case carries over to show that we can take $t_{i,n} = iT/n$
for the purposes of our calculation. The computation is then tedious
but straightforward, and carried out similarly to that for the
quadratic variation in the proof of Theorem 2 in \cite{twoscales}.
Theorem \ref{thm:discrerror} is thus proved.

\subsection{Proof of Corollary \ref{cor:discrerror}.} First of all,
note that since $M_n^3/n \rightarrow \infty$,
$\sum_{i=1}^{M_n} a_{n,i} = o ( - (n/M_n)^{1/2} )$. In lieu of equation
(\ref{eq:discrerror2}), it is therefore enough to prove
\begin{equation}
\label{eq:discrerror3}
(n/M_n)^{1/2} 
\sum_{i=1}^{M_n} a_{n,i} \left ( [X, X]^{(n,i)} - \XX \right ) \rightarrow
\eta_h Z
\end{equation}
Also, as in the proof of Theorem \ref{thm:discrerror}, our
assumptions imply that we can take $\mu_t = 0$ identically without
loss of generality.

Since there are asymptotically
infinitely
many $[X, X]^{(n,i)}$'s involved in equation (\ref{eq:discrerror2}),
we have to approximate with a finite number of these. To this end,
let $\delta > 0$ be an arbitrary number ($\delta < 1$). Let 
$\alpha = 1 - \delta / \sqrt{2}$. Let $L$ be an integer sufficiently
large that $2 \alpha^{L-1} \le \delta^2$. For $I = 1, ..., L$, let
$\tilde \kappa_I = \alpha^{L-I}$, and $\tilde \kappa_{0} = 0$.
For $i = 1, ..., M_n$, define $I_{i,n}$ to be the value $I , 1 \le I \le L$
for which $i/M_n \in (\tilde \kappa_{I-1} , \tilde\kappa_I ]$.
Then note that, if $|| U || = (EU^2)^{1/2}$,
\begin{equation}
\label{eq:discrerror4}
(n/M_n)^{1/2}
|| \sum_{i=1}^{M_n} a_{n,i} \left ( [X, X]^{(n,i)} - 
[X, X]^{(n,I_{i,n})}
\right ) ||
\le (n/M_n)^{1/2} \sum_{i=1}^{M_n} | a_{n,i} | 
\times \max_{1 \le i \le n} || [X, X]^{(n,i)} -
[X, X]^{(n,I_{i,n})} || .
\end{equation}
Now let $i_n$ be the value $i, 1 \le i \le M_n$ which maximizes
$|| [X, X]^{(n,i)} - [X, X]^{(n,I_{i,n})} ||$ for given $n$, and
let $I_n = I_{i_n,n}$. 

For the moment, let $N$ be an unbounded set of positive integers
so that $(i_n/M_n,I_n/M_n)_{n \in N}$ converges.
Call the limit $(\kappa_1,\kappa_2)$.
By the proof of Theorem \ref{thm:discrerror}
and of Theorem 2 in \cite{twoscales}, 
$(n/M_n) ( [X, X]^{(n,i_n)} - [X, X]^{(n,I_{n})} )^2$ is uniformly
integrable. By the statement of Theorem \ref{thm:discrerror}, it then 
follows that, as $n \rightarrow \infty$ through $N$
\begin{align}
(n/M_n) E ( [X, X]^{(n,i_n)} - [X, X]^{(n,I_{n})} )^2 &\rightarrow
E\eta^2 ( \Gamma_{2,2} + \Gamma_{1,1} - 2 \Gamma_{1,2} ) \notag \\
&= E\eta^2 2 T \kappa_2 \left ( 1 - \frac{\kappa_1}{\kappa_2} \right )^2
\notag \\
&\le E\eta^2 T \delta^2 
\label{eq:smdelta1}
\end{align}
by construction. Since every subsequence has a subsequence for which
$(i_n/M_n,I_n/M_n)$, it follows from equation (\ref{eq:discrerror4}) that
\begin{equation}
\label{eq:discrerror5}
\limsup_{n \rightarrow \infty} (n/M_n)^{1/2}
|| \sum_{i=1}^{M_n} a_{n,i} \left ( [X, X]^{(n,i)} -
[X, X]^{(n,I_{i,n})}
\right ) ||
\le \delta ( E\eta^2 T )^{1/2} \max_{0 \le x \le 1} | x h(x) |.
\end{equation}
The result of Corollary \ref{cor:discrerror} thus follows by
computing the limit of
\begin{equation}
\label{eq:discrerror6}
(n/M_n)^{1/2}
\sum_{i=1}^{M_n} a_{n,i} \left ( [X, X]^{(n,I_{i,n})} - \XX \right ) ,
\end{equation}
and then letting $\delta \rightarrow 0$.

\subsection{Proof of Theorem \ref{thm:MSRV-clt}.} The remainder term
$R_n$ from equation (\ref{eq:approx2}) can be written
$R_n = R_{n,1} + R_{n,2}$, where
\begin{equation}
R_{n,1} = \sum_{j=1}^{M_n} a_{n,j} \frac{1}{j} \left (\sum_{i=0}^{j-1}
\epsilon_{t_i}^2 + \sum_{i=n-j+1}^{n} \epsilon_{t_i}^2 \right )
-2 E \epsilon^2
\mbox{ and }
R_{n,2} = 2 \sum_{i=1}^{M_n} a_{n,i} [X, \epsilon]^{(i)}
\label{eq:R1and2}
\end{equation}
We shall show that $M_n^{1/2}R_n$ converges in law, 
conditionally on ${\cal X}$,
to a normal distribution with variance
\begin{equation}\label{eq:rem1}
4Var(\epsilon^2) \int_0^1 \int_0^y
x h(x){h(y)} dx dy  + 8 \XX Var(\epsilon) \int_0^1 \int_0^1 h(x)h(y) 
\min(x, y) dx dy ,
\end{equation}
and also that, conditionally
on ${\cal X}$, $R_n/M_n^{1/2}$ is asymptotically independent of
$(M_n^{3}/n)^{1/2}\zeta_n$ in Corollary \ref{cor:noise} in Section
\ref{sec:class}. Thus, in view of the results on the
the pure noise and discretization terms in
Corollaries \ref{cor:noise} and \ref{cor:discrerror},
Theorem \ref{thm:MSRV-clt} will then be shown. 

To show this, we show in the following that 
$M_n^{1/2}R_{n,1}$ and $M_n^{1/2}R_{n,2}$ are asymptotically normal
given ${\cal X}$, with mean zero and variances given by 
(\ref{eq:rem2}) and (\ref{eq:crossterm}), respectively. We then
discuss the {\it joint} distribution of $(M_n^{3}/n)^{1/2}\zeta_n$,
$M_n^{1/2}R_{n,1}$ and $M_n^{1/2}R_{n,2}$.

{\it Asymptotic normality of $R_{n,1}$.} Once $M_n < n/2$,
Write
\begin{equation}
R_{n,1} = \sum_{i=0}^{M_n-1} \epsilon_{t_i}^2 \sum_{j=i+1}^{M_n}
\frac{a_{n,j}}{j} + \sum_{i=0}^{M_n-1} \epsilon_{t_{n-i}}^2 \sum_{j=i+1}^{M_n}
\frac{a_{n,j}}{j} - 2E\epsilon^2.
\end{equation}
Hence,
\begin{eqnarray}
Var (M_n^{1/2}R_{n,1})
& =& 
2M_n Var(\epsilon^2) \sum_{i=0}^{M-1} ( \sum_{j=i+1}^{M} \frac{a_j}{j} )^2
\notag \\
&=& 2Var(\epsilon^2) \int_0^1 ( \int_x^1 h(y) dy)^2 dx + o(1) \notag \\
&=& 4Var(\epsilon^2) \int_0^1 \int_0^y x h(x)h(y) dx dy +o(1),
\label{eq:rem2}
\end{eqnarray}
while under Theorem \ref{th:hft},
\begin{eqnarray}
E \left [ \sum_{j=1}^{M_n} a_j \frac{1}{j} (\sum_{i=0}^{j-1}
\epsilon_{t_i}^2 + \sum_{i=n-j+1}^{n} \epsilon_{t_i}^2 ) \right ]=2
E\epsilon^2 (1 + o(M_n^{-1/2})).
\label{eq:rem.mean}
\end{eqnarray}
Since the Lindeberg condition is also obviously satisfied,
I obtain that $M_n^{1/2}R_{n,1}$ converges in law (conditionally on
${\cal X}$) to a normal distribution with mean zero and variance given
by equation (\ref{eq:rem2}).

{\it Asymptotic normality of the ``cross term'' $R_{n,2}$.}
As in the proof of Theorem \ref{thm:discrerror}, 
we proceed, without loss of generality, as if $X$ were a martingale.
Es in the proof of Theorem \ref{th:noise}, we shall show that
$M_n^{1/2}R_{n,2}$ satisfies the conditions of the version
of the Martingale Central Limit Theorem which is stated in
Corollary 3.1 (p. 58-59) of \cite{hallheyde}, and calculate
the asymptotic variance. As in the earlier proof,
we shall take, in the notation of
\cite{hallheyde}, ${\cal F}_{n, j}$ to be the smallest $\sigma$-field
making $\epsilon_{t_i}$, $i= 1, ... , j$, and the whole $X_t$ process,
measurable.

Note that, from (\ref{eq:revolK1}),
\begin{eqnarray*}
[X,\epsilon]^{(n,K)}
&=&\frac{1}{K}\sum_{i=0}^n  b_{n,i}^{(K)}\epsilon_{t_i},
\end{eqnarray*}
where
\[
b_{n,i}^{(K)}=\left\{
\begin{array}
[c]{ll}
-(X_{t_{n,i+K}} -X_{n,t_i}) &  \text{ \ if \ }  i = 0, \cdots, K-1\\
(X_{t_{n,i}} - X_{t_{n,i-K}}) - (X_{t_{n,i+K}} -X_{t_{n,i}})  &  \text{ \ if \ }
i = K, \cdots, n-K\\
(X_{t_{n,i}} - X_{t_{n,i-K}})  & \text{ \ if \ } i=n-K+1, \cdots, n
\end{array}
\right.
\]
Thus, from (\ref{eq:R1and2}), one obtains
\begin{align}
M_n^{1/2} R_{n,2} 
&= M_n^{1/2} \sum_{i=1}^n \epsilon_{t_i} \sum_{j=1}^{M_n} \frac{a_{n,j}}{j} 
b_{n,i}^{(j)} .
\end{align}
Obviously, $M_n^{1/2} R_{n,2}$ is the end point of a zero mean martingale 
relative to the filtration $({\cal F}_{n, j})$. The conditional variance
process (in Corollary 3.1 in \cite{hallheyde} is given by 
(we use $j\wedge k = \min(j,k)$)
\begin{align}
M_n E(\epsilon^2) \sum_{i=1}^n 
\left ( \sum_{j=1}^{M_n} \frac{a_{n,j}}{j} b_{n,i}^{(j)} \right )^2 &= 
M_n Var(\epsilon) \sum_{i=1}^n
\sum_{j=1}^{M_n} \sum_{k=1}^{M_n} \frac{a_{n,j}}{j} 
\frac{a_{n,k}}{k} b_{n,i}^{(j)} b_{n,i}^{(k)}
\notag \\
&=M_n Var(\epsilon) \sum_{i=1}^n
\sum_{j=1}^{M_n} \sum_{k=1}^{M_n} \frac{a_{n,j}}{j}
\frac{a_{n,k}}{k} (b_{n,i}^{(j\wedge k)} )^2 +o_p(1)
\notag \\
&= 2 M_n Var(\epsilon) \sum_{j=1}^{M_n} \sum_{k=1}^{M_n} \frac{a_{n,j}}{j}
\frac{a_{n,k}}{k} (j\wedge k) [X,X]^{(j\wedge k)} +o_p(1) \notag \\
& =  2 \int_0^1 \int_0^1 h(x)h(y) (x\wedge y) dx dy \XX Var(\epsilon) 
+ o_p(1 ), 
\label{eq:crossterm}
\end{align}
where remainder terms are taken care of as in the proof of Theorem
\ref{thm:discrerror}.

By similar methods, the Lindeberg condition is satisfied (cf. the
discussion in the proof of Theorem \ref{th:noise}). By Corollary 3.1
(p. 58-59) in \cite{hallheyde}
it follows that $M^{1/2}_n R_n$
is asymptotically normal (conditionally on ${\cal X}$), 
with mean zero and variance given by
(\ref{eq:crossterm}). This is what we needed to show. 

{\it The joint distribution of $(M_n^{3}/n)^{1/2}\zeta_n$,
$M_n^{1/2}R_{n,1}$ and $M_n^{1/2}R_{n,2}$.} First of all, note that
for all three quantities, we have satisfied the conditions of 
Corollary 3.1 (p. 58-59) of \cite{hallheyde}. This is with the exception of 
(their equation) (3.21), where we have instead used the Remarks following
their corollary (and thus the convergence is conditional on ${\cal X}$
as opposed to stable with respect to the $\sigma$-field generated
by both ${\cal X}$ and the $\epsilon_{t_i}$).

In terms of joint distribution, note first that the sum of conditional
covariances (for each two of the three quantities 
$(M_n^{3}/n)^{1/2}\zeta_n$, $M_n^{1/2}R_{n,1}$ and $M_n^{1/2}R_{n,2}$
converge to zero, by the same methods as above. In view of how
Hall and Heyde's corollary implies their Theorem 3.2 (p. 58), 
the Cram{\'e}r-Wold device now implies the joint normality
of $(M_n^{3}/n)^{1/2}\zeta_n$, $M_n^{1/2}R_{n,1}$ and $M_n^{1/2}R_{n,2}$,
and also that they are asymptotically independent.
Theorem \ref{thm:MSRV-clt} is then proved.

\bibliographystyle{asa}
\bibliography{mainbib.bib}

\end{document}